\documentclass[preprint,12pt]{amsart}
\usepackage{amscd,amssymb,amsthm,latexsym,stmaryrd,xypic,amsmath}
\usepackage{mathrsfs,yhmath}

\newcommand{\M}{\mathrm{M}}

\newcommand{\Tr}{\mathrm{Tr}}
\newcommand{\Mat}{\mathrm{Mat}}

\def\ov#1{{\overline{#1}}}

\newcommand{\ds}{\displaystyle}

\renewcommand{\to}{\longrightarrow}

\newcommand{\funcp}[4]{\begin{aligned}\newline #1&\longrightarrow #2 \cr\newline #3 &\longmapsto #4.\end{aligned}}

\newcommand{\nfuncp}[5]{#1\colon\begin{aligned}\newline #2&\longrightarrow #3 \cr\newline #4 &\longmapsto #5.\end{aligned}}
\newcommand{\nfuncv}[5]{#1\colon\begin{aligned}\newline #2&\longrightarrow #3 \cr\newline #4 &\longmapsto #5,\end{aligned}}

\let\cdotorg\cdot \def\cdot{{\cdotorg}}

\xyoption{all}

\numberwithin{equation}{section}
\theoremstyle{plain} %
\newtheorem{thm}[equation]{Theorem}
\newtheorem{prop}[equation]{Proposition}
\newtheorem{coro}[equation]{Corollary}
\newtheorem{lem}[equation]{Lemma}
\theoremstyle{definition} %
\newtheorem{ex}[equation]{Example}

\newtheorem{rem}[equation]{Remark}

\renewenvironment{proof}{{\it Proof.}}{
  \qed\endtrivlist
  \par}

\usepackage{hyperref}
\hypersetup{
hyperindex=true, 
colorlinks=true, 
linkcolor=blue,
linktocpage=true,           
filecolor=red,
citecolor=red,
breaklinks=false, 
urlcolor= blue, 
bookmarks=true, 
bookmarksopen=false,
}

\title[Symmetric matrices in characteristic two]{Minimal and characteristic polynomials of symmetric matrices in characteristic two}

\author{Gr\'{e}gory Berhuy}

\begin{document}

\maketitle

 \address{Institut Fourier, Universit\'e Grenoble Alpes, 100, rue des maths 38610 Gières, France}
 
{\bf Email : } gregory.berhuy@univ-grenoble-alpes.fr

\begin{abstract}
Let $k$ be a field of characteristic two. We prove that a monic polynomial $f\in k[X]$ of degree $n\geq 1$ is the minimal/characteristic polynomial of a symmetric matrix with entries in $k$ if and only if it is not the product of pairwise distinct inseparable irreducible polynomials. In this case, we prove that $f$ is the minimal polynomial of a symmetric matrix of size $n$. We also prove that any element $\alpha\in k_{alg}$ of degree $n\geq 1$ is the eigenvalue of a symmetric matrix of size $n$ or $n+1$, the first case happening if and only if the minimal polynomial of $\alpha$ is separable. 
\end{abstract}



{\it Keywords : }Symmetric matrices, \  Minimal polyomial, \ Characteristic polynomial, \ Eigenvalues, \ Symmetric bilinear forms, \  Transfer

{\it 2020 MSC Codes: } 11C20, \ 15A15, \ 15A18, \ 11E39


\tableofcontents

\bigskip

{\bf Notation. }The following notation will be used throughout this paper.

\begin{itemize}

\item If $k$ is a field, $k_{alg}$ will denote a fixed algebraic closure of $k$, and $k_s$ will denote the separable closure of $k$ in $k_{alg}$.

\item Let $V$ be a finite dimensional $k$-vector space, and let $\mathcal{B}$ be a basis of $V$. If $u:V\to V$ is an endomorphism of $V$,  its matrix representation with respect to the basis $\mathcal{B}$  will be denoted by $\Mat(u;\mathcal{B})$.

Similarly, the Gram matrix of a symmetric $k$-bilinear form $b:V\times V\to k$ with respect to the basis $\mathcal{B}$ will be denoted by $\Mat(b;\mathcal{B})$.

\item For $n\geq 1$, $\M_n(k)$ will denote the ring of $n\times n$ matrices with entries in $k$. If $M\in \M_n(k)$, the minimal polynomial and the characteristic polynomial of $M$ will be denoted respectively by $\mu_M$ and $\chi_M.$
\end{itemize}

\section{Introduction}
The problem of determining the minimal and characteristic polynomials of symmetric matrices with entries in a field  has a long history, that we briefly sketch. The first major result on this question is due to Krakowski. If $k$ is a formally real field, we say that $f\in k[X]$ is totally real 
if it splits over every real closure of $k$. Krakowski proved (\cite{Kra}) that if $k$ is a formally real field, a non constant monic polynomial $f\in k[X]$ is the minimal polynomial of a symmetric matrix with entries in $k$ if and only if it is separable and totally real. He also proved that if $k$ is not formally real and ${\rm char}(k)\neq 2$, any non constant monic polynomial $f\in k[X]$  is  the minimal polynomial of a symmetric matrix with entries in $k$. In particular, any totally real algebraic number (resp. any algebraic number) $\alpha\in k_{alg}$ is the eigenvalue of a symmetric matrix with entries in $k$ if $k$ is formally real (resp. if $k$ is not formally real of characteristic different from two).
However, the size of the symmetric matrices constructed by Krakowski was huge compared to the degree of $f$, and the question of finding the minimal size of such symmetric matrices was still open. This question was solved by Bender. Concerning the eigenvalues of symmetric matrices, he proved (\cite{Ben1}) that, if $k$ is a number field, any totally real number of degree $n$ is the eigenvalue of a symmetric matrix of $\M_n(k)$ or $\M_{n+1}(k).$
When ${\rm char}(k)\neq 2$, given a monic polynomial $f\in k[X],$ (where $f$ is totally real if $k$ is formally real), Bender (\cite{Ben2}) computed the smallest integer $r$ such that there exists a symmetrix matrix with entries in $k$ satisfying $\mu_M=f$ and $\chi_M=f^r$.
In particular, if $-1$ is a square in $k$, Bender's results imply that any monic polynomial $f\in k[X]$ of degree $n\geq 1$ is the minimal/characteristic polynomial of a symmetric matrix of $\M_n(k)$.

When ${\rm char}(k)=2$, much less is known. In \cite{Ben2}, Bender stated without proof the following result: if $k$ has characteristic two, and if $f\in k[X]$ is a monic polynomial of degree $n$ which has at least one separable irreducible divisor or such that all the valuations corresponding to inseparable irreducible divisors are even, then $f$ is the minimal polynomial of a symmetric matrix of $\M_n(k)$.

Finally, results of Bass, Estes and Guralnick (see \cite{BEG}, (9.4)) imply that, if  $k$ has characteristic two, any element $\alpha\in k_{alg}$ of degree $n$ is the eigenvalue of a symmetric matrix with entries in $k$ of size at most $2n+1$.

In this paper, we will determine the minimal and characteristic polynomials of symmetric matrices with entries in a field of characteristic two. More precisely, we will prove the following results.

\begin{thm}
Let $k$ be a field of characteristic two, and let $f\in k[X]$ be a monic polynomial of degree $n\geq 1$. Then  $f$ is the minimal polynomial of a symmetric matrix with entries in $k$ if and only if $f$ is not the product of pairwise distinct monic irreducible inseparable polynomials.

In this case, $f$ is the minimal polynomial of a symmetric matrix of $\M_n(k)$.

In particular, if $k$ is perfect, any monic polynomial $f\in k[X]$ of degree $n\geq 1$  is the minimal polynomial of a symmetric  matrix of $\M_n(k)$.
\end{thm}

\begin{coro}
Let $k$ be a field of characteristic two, and let $f\in k[X]$ be a monic polynomial of degree $n\geq 1$. Then $f$ is the characteristic polynomial of a symmetric matrix of $\M_n(k)$ if and only if it is not the product of pairwise distinct monic irreducible inseparable polynomials.

In particular, if $k$ is perfect, any monic  polynomial $f\in k[X]$ of degree $n\geq 1$ is the characteristic polynomial of a symmetric matrix of $\M_n(k)$.
\end{coro}

\begin{thm}
Let $k$ be a field of characteristic two, and let $\alpha\in k_{alg}$ be an algebraic element of degree $n$, with minimal polynomial $f.$ Then :

\begin{enumerate}
 \item if $f$ is separable, $\alpha$ is the eigenvalue of a symmetric matrix of $\M_n(k)$;
 
 \item if $f$ is inseparable, $\alpha$ is the eigenvalue of a symmetric matrix of $\M_{n+1}(k)$, but not of any symmetric matrix of $\M_n(k)$.
\end{enumerate}

In particular, if $k$ is perfect, any algebraic element of degree $n$ is the eigenvalue of a symmetric  matrix of $\M_n(k)$.
\end{thm}

The proofs rely on a technique introduced by Bender, which consists in constructing $k$-linear forms $k[X]/(f)\to k$ such that the corresponding transfer is isomorphic to the unit form (see Lemma \ref{kxftransfer}).

\section{Bilinear forms and symmetric matrices}\label{secbil}

We first recall some well-known definitions on bilinear forms.

The {\it unit form of rank $n$} is the symmetric $k$-bilinear form $$\funcp{k^n\times k^n}{k}{(x,y)}{x^t y}$$

A symmetric $k$-bilinear form $b:V\times V\to k$ on a finite dimensional $k$-vector space will then be isomorphic to the unit form (of rank $\dim_k(V)$) if and only if $V$ has an orthonormal basis with respect to $b$.

If $E$ is a commutative $k$-algebra and $s:E\to k$ is a $k$-linear form, the {\it transfer} associated to $s$ is the symmetric $k$-bilinear form $$\nfuncp{s_*(\langle 1\rangle)}{E\times E}{k}{(x,y)}{s(xy)}$$ 

We now prove a slight variation of a lemma of Bender, which relates the theory of symmetric matrices and the notion of transfer. 

Before stating it, recall that if $E$ is a finite dimensional commutative $k$-algebra, the {\it minimal polynomial} of $x\in E$ is the unique monic generator $\mu_x$ of the ideal $I_x$ of $k[X]$ defined by  $$I_x=\{P\in k[X]\mid P(x)=0\}.$$

\begin{lem}\label{kxftransfer}
Let $k$ be an arbitrary field, and let $E$ be a commutative $k$-algebra of dimension $n\geq 1.$ Assume that there exists a $k$-linear form $s:E\to k$ such
that $s_*(\langle 1\rangle)$ is isomorphic to the unit form, and let $\mathcal{B}$ be an orthonormal basis of $E$. Then for all $x\in E$, the matrix 
$M_x={\rm Mat}(\ell_x;\mathcal{B})$ is a symmetric matrix of $\M_n(k)$ such that $\mu_{M_x} =\mu_x,$ where $\ell_x:E\to E$ denotes the endomorphism of left multiplication by $x$. 
\end{lem}

\begin{proof}
Assume that there exists a $k$-linear form $s:E\to k$ such
that $s_*(\langle 1\rangle)$ is isomorphic to the unit form, and let $\mathcal{B}$ be an orthonormal basis of $E$ with respect to $s_*(\langle 1\rangle)$.
 Then for all $x,x_1,x_2\in E$, we have $$s(\ell_x(x_1)x_2)=s(x x_1x_2)=s(x_1  xx_2)=s(x_1\ell_x(x_2)).$$  Hence $\ell_x$ is self-adjoint with respect to $s_*(\langle 1\rangle)$. It follows that $M_x=\Mat(\ell_x; \mathcal{B})$ is symmetric. 
Now for all $P\in k[X]$, we have $$P(\ell_x)=0\iff \ell_{P(x)}=0\iff P(x)=0\iff \mu_x\mid P.$$ This implies that the minimal polynomial of $M_x$ is $\mu_x$.
\end{proof}

For the rest of the paper, we will assume that $k$ is a field of characteristic two.  Recall that a $k$-bilinear form $b:V\times V\to k$ is {\it alternating} if $b(x,x)=0$ for all $x\in V$. Alternating forms are necessarily symmetric. If $V$ is finite dimensional, non-degenerate alternating $k$-bilinear forms are {\it hyperbolic}, that is, isomorphic to an orthogonal sum of hyperbolic planes $\mathbb{H}$, where $$\nfuncp{\mathbb{H}}{k^2\times k^2}{k}{(x,y)}{x_1y_2+x_2y_1}$$

The following lemma gives a nice characterization of the unit form in characteristic two.

\begin{lem}\label{unit}
A nonzero symmetric bilinear form $b:V\times V\to k$ is isomorphic to the unit form if and only it is non-degenerate, non-alternating, and $b(x,x)$ is a square for all $x\in V$.

In particular, if $m\geq 1$ and $h$ is a hyperbolic form, $m\times \langle 1 \rangle \perp h$ is isomorphic to the unit form.
\end{lem}

\begin{proof}
Since $k$ has characteristic two, any sum of squares is a square, and the unit form satisfies the conditions of the lemma. Conversely, assume that  $b:V\times V\to k$ is a nonzero symmetric bilinear form satisfying the conditions of the lemma. Then there exists $e_1\in V$ such that $b(e_1,e_1)\neq 0$, and by assumption $b(e_1,e_1)=\lambda^2$ for some $\lambda\in k^\times$. Replacing $e_1$ by $\lambda^{-1}e_1$ if necessary, one may assume that $b(e_1,e_1)=1$.
The restriction of $b$ to $F=k e_1$ being non-degenerate, we have $V=F\oplus F^\perp$. Now, the restriction of $b$ to $F^\perp$ is also non-degenerate (if $(e_2,\ldots,e_m)$ is a basis of $F^\perp$, then the matrix of $b$ in the basis $(e_1,\ldots,e_n)$ is block-diagonal, with a $1$ in the left upper corner).
If this restriction is non-alternating, one may apply induction on $\dim_k(V)$ (the case $\dim_k(V)=1$ being trivial).
If it is alternating, then it is hyperbolic. Hence, it is enough to prove that the matrix $B=\begin{pmatrix}
1 & 0 & 0 \cr 0 & 0 & 1\cr 0 & 1 & 0
\end{pmatrix}$ is congruent to the identity matrix, and apply induction. But one may check that $P^tBP=I_3$, where $P=\begin{pmatrix}
1 & 1 & 1\cr 1 & 0 & 1 \cr 0 & 1 & 1
\end{pmatrix}.$ The last part is clear.
\end{proof}

\begin{rem}
The previous lemma is a direct consequence of the following well-known fact, which may be proved using a slight modification of the arguments above: any non-alternating symmetric bilinear space $(V,b)$ has an orthogonal basis.
\end{rem}

\begin{rem}\label{squares}
Let $b:V\times V\to k$ be a symmetric bilinear form, and let $\mathcal{B}=(e_1,\ldots,e_n)$ be a basis of $V.$ If $x=\ds\sum_{i=1}^n x_i\cdot e_i$, since $k$ has characteristic two, we have $b(x,x)=\ds\sum_{i=1}^n x_i^2 b(e_i,e_i)$. Using the fact that a sum of squares is a square, we see that the following properties are equivalent:

\begin{enumerate}
\item $b(x,x)$ is a square for all $x\in V$;

\item for any basis $\mathcal{B}=(e_1,\ldots,e_n)$, $b(e_i,e_i)$ is a square for all $1\leq i\leq n$;

\item there exists a basis  $\mathcal{B}=(e_1,\ldots,e_n)$ such that $b(e_i,e_i)$ is a square for all $1\leq i\leq n$.
\end{enumerate}

The same argument shows that the following properties are equivalent:

\begin{enumerate}
\item $b$ is alternating;

\item for any basis $\mathcal{B}=(e_1,\ldots,e_n)$, $b(e_i,e_i)=0$  for all $1\leq i\leq n$;

\item there exists a basis  $\mathcal{B}=(e_1,\ldots,e_n)$ such that $b(e_i,e_i)=0$ for all $1\leq i\leq n$.
\end{enumerate}

This remark will be useful to prove that a bilinear form is isomorphic to the unit form. 
\end{rem}

For practical computations, we now explain a modified version of Gauss reduction  which computes an orthonormal basis of a bilinear space $(V,b)$ satisfying the conditions of Lemma \ref{unit}.

If $\varphi, \psi\in V^*$, we will denote by $\varphi\bullet \psi$ the bilinear form $$\nfuncp{\varphi\bullet\psi}{V\times V}{k}{(x,y)}{\varphi(x)\psi(y)}$$

The proposed algorithm will rely on the following identities : if $\varphi_1,\varphi_2,\varphi_3\in V^*$, we have 

\begin{enumerate}
\item[(i)] $\varphi_1\bullet\varphi_1+\varphi_1\bullet \varphi_2+\varphi_2\bullet\varphi_1=(\varphi_1+\varphi_2)\bullet(\varphi_1+\varphi_2)+\varphi_2\bullet\varphi_2$

\item [(ii)] $\varphi_1\bullet\varphi_1+\varphi_2\bullet \varphi_3+\varphi_3\bullet\varphi_2=(\varphi_1+\varphi_2)\bullet(\varphi_1+\varphi_2)+(\varphi_1+\varphi_3)\bullet(\varphi_1+\varphi_3)+(\varphi_1+\varphi_2+\varphi_3)\bullet(\varphi_1+\varphi_2+\varphi_3)$
\end{enumerate}

Let $b:V\times V\to k$ be a non-alternating symmetric bilinear form such that $b(x,x)$ is a square for all $x\in V$. We do not assume here that $b$ is non-degenerate. Our first goal is to write $b=\ds\sum_{i=1}^r\varphi_i\bullet\varphi_i$, where $\varphi_1,\ldots,\varphi_r\in V^*$ are $k$-linearly independent linear forms on $V$.
 
We start with the Gram matrix $M=(a_{ij})$ of $b$ in a fixed basis $\mathcal{E}=(e_1,\ldots,e_n)$ of $V$. 

If $x=\ds\sum_{i=1}^n x_i\cdot e_i, y=\ds\sum_{j=1}^n y_j \cdot e_j$, we then have $$b(x,y)=\ds\sum_{i=1}^n a_{ii}x_iy_i+\sum_{i<j} a_{ij}(x_iy_j+x_j y_i).$$

{\bf Step 1. }Since $b$ is non-alternating, it contains a term of the form $u^2x_iy_i,$ where $u\in k^\times$. For sake of simplicity, assume that $i=1$. Collecting all terms containing $x_1$ or $y_1$, write $$b(x,y)=(u x_1)(uy_1)+(ux_1)\varphi(y)+(uy_1)\varphi(x)+ \mbox{ the remaining terms},$$ where $\varphi$ is a linear combination of $e_2^*,\ldots,e_n^*$. If we set $$\nfuncv{\varphi_1}{V}{k}{x}{ux_1+\varphi(x)}$$
using identity (i), we may write $b=\varphi_1\bullet \varphi_1+b_1,$ where $b_1$ is a symmetric bilinear form on $V$, whose expression only depends on $x_2,y_2,\ldots,x_n,y_n$. 

Note that $b_1(x,x)$ is a square for all $x\in V$, since $b_1(x,x)=\varphi_1(x)^2+ b(x,x)$ and a sum of squares is a square.

{\bf Step 2. }If $b_1=0$, we are done. Otherwise, if $b_1$ is non-alternating, we repeat Step 1 with $b_1$. If $b_1$ is alternating, it contains a nonzero term of the form $a_{ij}(x_i y_j+x_j y_i).$ Say for example that $(i,j)=(2,3)$. Collecting terms containing $x_2,x_3,y_2$ or $y_3$, we may write 
$$b_1(x,y)=a_{23}(x_2y_3+x_3y_2)+\psi(x)y_3+\psi(y)x_3+\psi'(x)y_2+\psi'(y)x_2+$$ 
$$\mbox{ the remaining terms in }x_i,y_i,i\geq 4$$
where $\psi,\psi'$ are linear combinations of $e_4^*,\ldots,e_n^*$.

Now, $a_{23}(x_2y_3+x_3y_2)+\psi(x)y_3+\psi(y)x_3+\psi'(x)y_2+\psi'(y)x_2$ is equal to 
$$(a_{23}x_2+ \psi(x) )(y_3+ a_{23}^{-1}\psi'(y)) + (a_{23} y_2+ \psi(y) )(x_3+ a_{23}^{-1}\psi'(x))+a_{23}^{-1}(\psi(x)\psi'(y)+\psi'(x)\psi(y)).$$

Hence, we may write $$b=\varphi_1\bullet\varphi_1+\psi_2\bullet\psi_3+\psi_3\bullet \psi_2+b_2,$$
where $b_2$ is a symmetric bilinear form, whose expression only depends on $x_4,y_4,\ldots,x_n,y_n.$ 
Using identity (ii), and replacing the former  $\varphi_1$ by $\varphi_1+\psi_2$,   we have  $b=\ds\sum_{j=1}^3\varphi_j\bullet \varphi_j+b_2$, for some $\varphi_1,\varphi_2,\varphi_3\in V^*$.

Once again, $b_2(x,x)=0$ for all $x\in V$. If $ b_2=0$, we are done. Otherwise, we repeat Step 1 or Step 2 with $b_2$.

Since the number of variables decrease by one or two at each step, we will end with a form $b_1$ or $b_2$ which is identically zero in at most $r$  iterations of Steps $1$ or $2$.
 
Using induction on $n$, it is not difficult to prove, as in the classical Gauss reduction algorithm, that the linear forms $\varphi_1,\ldots,\varphi_r$ obtained at the end of this procedure are linearly independent.

{\bf Step 3. }We now have $b=\ds\sum_{j=1}^r\varphi_j\bullet\varphi_j,$ where $\varphi_1,\ldots,\varphi_r\in V^*$ are linearly independent.

If $\varphi_i=\ds\sum_{j=1}^n \varphi_{ij}e^*_j$, the matrix $U=(\varphi_{ij})\in\M_{r\times n}(k)$  has rank $r$. We may then add $n-r$ rows to $U$ in order to obtain a  matrix $Q\in{\rm GL}_n(k)$. Note that this step is not necessary if $b$ is non-degenerate, that is, if $r=n$.

The columns of $P=Q^{-1}$ then represent the  coordinate vectors of the elements of an orthogonal basis $\mathcal{B}$ in the basis $\mathcal{E},$ for which the corresponding Gram matrix is $$\Mat(b;\mathcal{B})=\begin{pmatrix}
I_r &0 \cr 0 & 0\end{pmatrix}.$$
In particular, if $b$ is non-degenerate, $\mathcal{B}$ will be an orthonormal basis of $V$ with respect to $b$.

\begin{rem}
The previous procedure may be slightly modified to produce an orthogonal basis of any non-alternating symmetric bilinear space $(V,b)$.
\end{rem}

We now go back to our original problem.

\begin{lem}\label{subunit}
Let $b:V\times V\to k$ be a nonzero symmetric bilinear form on a finite dimensional $k$-vector space $V.$ 

Assume that $(V,b)\simeq (V_1,b_1)\perp \cdots \perp (V_r, b_r)$, where $V_1,\ldots,V_r$ are nonzero vector spaces and $b_1,\ldots,b_r$ are symmetric bilinear forms. Then $b$ is isomorphic to the unit form if and only if the following properties hold:

\begin{enumerate}
\item for all $1\leq i\leq r$, $b_i$ is either isomorphic to the unit form or hyperbolic;

\item there exists at least one $1\leq j\leq r$ such that $b_j$ is isomorphic to the unit form.
\end{enumerate} 
\end{lem}

\begin{proof}
If $b_1,\ldots,b_r$ satisfy the conditions of the lemma, then $b$ is isomorphic to the unit form by Lemma \ref{unit}, since $b\simeq b_1\perp\cdots\perp b_r$ is isomorphic to the  orthogonal sum of a unit form and of an hyperbolic form. Conversely, assume that $b$ is isomorphic to the unit form. Then $b_1,\ldots,b_r$ are all non-degenerate and nonzero, since $b$ is non-degenerate and each $V_i$ is nonzero. If $b_i$ is alternating, $b_i$ is hyperbolic. Assume now that $b_i$ is non-alternating. Since $b(x,x)$ is a square for all $x\in V$, the same property holds for $(V_i, b_i)$. In this case, $b_i$ is isomorphic to the unit form by Lemma \ref{unit} again. 
Finally, if $b_1,\ldots,b_r$ are all hyperbolic, so is $b$. In this case, $b$ is alternating, which contradicts the fact that $b$ is isomorphic to the unit form. This concludes the proof.
\end{proof}

The previous lemma, applied to transfers, yields the following result.

\begin{lem}\label{sumtrans}
Let $E_1,\ldots,E_r$ be $r$ nonzero finite dimensional commutative $k$-algebras. Then there exists a $k$-linear form $s:E_1\times \cdots \times E_r\to k$ such that $s_*(\langle 1\rangle)$ is isomorphic to the unit form if and only if there exist $r$ $k$-linear forms $s_i:E_i\to k$ such that the following properties hold:

\begin{enumerate}
\item for all $1\leq i\leq r$, $(s_i)_*(\langle 1 \rangle)$ is either isomorphic to the unit form or hyperbolic;

\item there exists at least one $1\leq j\leq r$ such that  $(s_j)_*(\langle 1 \rangle)$ is isomorphic to the unit form.
\end{enumerate} 
\end{lem}

\begin{proof}
Any $k$-linear form $s:E_1\times \cdots \times E_r\to k$ may be written in a unique way as $$\nfuncv{s=s_1\oplus\cdots\oplus s_r}{E_1\times \cdots\times E_r}{k}{(x_1,\ldots,x_r)}{s_1(x_1)+\cdots+s_r(x_r)}$$
where $s_i:E_i\to k$ is a $k$-linear form. For such a decomposition, we clearly have $$s_*(\langle 1\rangle)\simeq (s_1)_*(\langle 1\rangle)\perp \cdots \perp (s_r)_*(\langle 1\rangle).$$ Now apply Lemma \ref{subunit} to conclude.
\end{proof}

From this lemma, we now derive some information on the structure of minimal polynomials of symmetric matrices.

\begin{prop}\label{prodrho}
Let $\rho_1,\ldots,\rho_s\in k[X]$ be pairwise distinct monic inseparable irreducible polynomials. Then $\rho_1\cdots\rho_s$ is not the minimal polynomial of a symmetric matrix with entries in $k$ of any size.
\end{prop}

\begin{proof}
Assume that there exists a symmetric matrix $M\in \M_n(k)$ such that $\mu_M=\rho_1\cdots\rho_s$. 

Write $\rho_i=\pi_i(X^{2^{k_i}}),$ where $k_i\geq 1$ and $\pi_i\in k[X]$ is an irreducible separable polynomial. Thus,  $\rho_i=\ds\prod_{j=1}^{d_i}(X^{2^{k_i}}-a_{ij})\in k_s[X],$ where $a_{i1},\ldots,a_{id_i}\in k_s$ are the roots of $\pi_i.$

Hence $$\mu_M=\ds\prod_{i=1}^s\prod_{j=1}^{d_i}(X^{2^{k_i}}-a_{ij})\in k_s[X],$$
and $k_s^n=\ds\bigoplus_{i,j}E_{ij}$, where $E_{ij}=\ker(M^{2^{k_i}}-a_{ij} I_n)$.

Let us prove that the subspaces $E_{ij}$ are mutually orthogonal with respect to the unit form on $k_s^n$.
Notice that, if $(i,j)\neq (\ell,m)$, the polynomials $P=X^{2^{k_i}}-a_{ij}$ and $Q=X^{2^{k_\ell}}-a_{\ell m}$ are coprime : if $i\neq \ell$, it follows from the fact that $\rho_i$ and $\rho_\ell$ are distinct irreducible polynomials (and therefore coprime), and if $i=\ell$ and $j\neq m$, this comes from the fact $a_{ij}\neq a_{im}$ since $\pi_i$ is separable, and from the relation $P-Q=a_{ij}-a_{im}$.

Thus, there exists $U,V\in k_s[X]$ such that $UP+VQ=1$. Now, for all $x\in E_{ij}=\ker(P(M))$ and all $y\in E_{\ell m}=\ker(Q(M))$, since $M$ is symmetric, we have 
$$\begin{array}{lllll}x^ty&=&(I_nx)^t y&=&(U(M)P(M)x+V(M)Q(M)x)^ty\cr
 &=&(V(M)Q(M)x)^ty&=&x^t(Q(M)V(M)y)\cr &=&x^t(V(M)Q(M)y)&=&0.\end{array}$$
Thus, the subspaces $E_{ij}$ are mutually orthogonal and therefore, the unit form $b_0$ on $k_s^n$ satisfies $b_0\simeq \perp_{i,j} b_{ij},$ where $b_{ij}$ is the restriction of $b_0$ to $E_{ij}$. By Lemma \ref{subunit}, there exist $i,j$ such that $b_{ij}$ is isomorphic to the unit form. In other words, $E_{ij}$ has an orthonormal basis for the restriction of $b_0$ to $E_{ij}$. Now $E_{ij}$ is stable by $M$, so $M$ induces an endomorphism $u_{ij}$ on $E_{ij}$, which is self-adjoint with respect to $b_{ij}$. Hence the matrix $M'$ of $u_{ij}$ in an orthonormal basis of $E_{ij}$ is a symmetric matrix with entries in $k_s$ satisfying $(M')^{2^{k_i}}=a_{ij}I_m,$ where $m$ is the size of $M'$. If $N=(M')^{2^{k_i-1}}$, then $N=(n_{ij})$ is symmetric and satisfies $N^2=a_{ij}I_m$. Comparing the coefficients in position $(1,1)$, we get $$a_{ij}=\sum_{\ell=1}^mn_{1\ell}n_{\ell 1}=\sum_{\ell=1}^mn_{1\ell}^2=(\sum_{\ell=1}^m n_{1\ell})^2.$$ Therefore, $a_{ij}$ is a square in $k_s$. 
Since the Galois group  of a splitting field of $\pi_i$ acts transitively on its roots, this implies  that all the roots of $\pi_i$ are squares in $k_s$. It follows that all the coefficients of $\rho_i=\pi_i(X^{2^{k_i}})$ are squares of $k_s$ (as a sum of products of squares).  Now, if $a\in k$ is a square of an element of $k_s$, it is already a square of an element of $k$. Indeed, if $\lambda\in k_s$ satisfies $\lambda^2=a$, then $\mu_{\lambda,k}\mid X^2-a=(X-\lambda)^2$. Since $\lambda$ is separable over $k$, the previous divisibility relation forces $\mu_{\lambda,k}=X-\lambda$, and $\lambda\in k.$ It follows that $\rho_i$ is a square in $k[X],$ contradicting its irreducibility. This concludes the proof.
\end{proof}

\section{Computation of various transfers}

In the sequel, if $f\in k[X]$ is a monic polynomial of degree $n\geq 1,$ $\alpha$ will denote the class of $X$ in $k[X]/(f)$ (unless specified otherwise), so that $(1,\alpha, \ldots, \alpha^{n-1})$  is a $k$-basis of $k[X]/(f)$.

\begin{lem}\label{even}
Assume that $f\in k[X^2]$ has degree $2m$, and let $s:k[X]/(f)\to k$ be the unique $k$-linear form such that $$s(1)=\cdots=s(\alpha^{2m-2})=0, s(\alpha^{2m-1})=1.$$
Then $s_*(\langle 1\rangle)$ is hyperbolic.
\end{lem}

\begin{proof}
It is enough to prove that $s_*(\langle 1\rangle)$ is non-degenerate and alternating.

Let us prove it is non-degenerate. Let $x=\ds\sum_{i=0}^{2m-1}\lambda_i\cdot \alpha^i\in \ker(s_*(\langle 1\rangle))$. 
Then $s(x 1)=s(x)=\lambda_{2m-1}=0$. Assume we proved that $\lambda_{2m-1}=\lambda_{2m-2}=\cdots=\lambda_{2m-1-j}=0$ for some $0\leq j\leq 2m-2$.
Then, $$0=s(x \alpha^{j+1})=s(\ds\sum_{i=0}^{2m-j-2}\lambda_i\cdot \alpha^{i+j+1}))=s(\ds\sum_{i=j+1}^{2m-1}\lambda_{i-j-1}\cdot \alpha^i))=\lambda_{2m-1-j-1}.$$
It follows by induction that all the $\lambda_i$'s are zero, that is $x=0$. Hence, $s_*(\langle 1\rangle)$ is non-degenerate.

We now prove that it is alternating. By Remark \ref{squares}, it is enough to show that $s(\alpha^{2i})=0$ for $0\leq i\leq 2m-1$. For, let $R$ be  the remainder of the long division of $X^{2i}=(X^2)^{i}$ by $f$. It is clear that $R\in k[X^2]$ (since we may just think of $X^2$ as the variable when performing the long division).

Consequently, evaluation at $\alpha$ shows that $\alpha^{2i}=R(\alpha)$ is a linear combination of $1,\alpha^2,\ldots,\alpha^{2m-2}$, and $s(\alpha^{2i})=0$.
This concludes the proof.
\end{proof}

\begin{lem}\label{transtr}
Let $L/k$ be a finite separable field extension, and let $\Tr_{L/k}:L\to k$ be the corresponding trace map. Then $(\Tr_{L/k})_*(\langle 1\rangle)$ is isomorphic to the unit form.
\end{lem}

\begin{proof}
If $X_L$ denotes the set of $k$-embeddings $L\to k_{alg}$, recall that for all $x\in L$, we have $\Tr_{L/k}(x)=\ds\sum_{\sigma\in X_L}\sigma(x).$
In particular, for all $x\in L$, we get $$\Tr_{L/k}(x^2)=\sum_{\sigma\in X_L}\sigma(x^2)=\sum_{\sigma\in X_L}\sigma(x)^2=\Tr_{L/K}(x)^2.$$
Since $L/k$ is separable, it is known that $\Tr_{L/k}$ is nonzero, and that the corresponding transfer is non-degenerate. Since the trace map is nonzero, the previous equality shows that the corresponding transfer is non-alternating. It also implies  that $\Tr_{L/k}(x^2)$ is a square for all $x\in L$. Lemma \ref{unit} then yields the desired result. 
\end{proof}

\begin{lem}\label{xam}
Assume that $f=(X-a)^m$, $m\geq 1$, and let $s:k[X]/(f)\to k$ be the unique $k$-linear form such that $$s(1)=s(\alpha-a)=\cdots=s((\alpha-a)^{m-1})=1.$$
Then $s_*(\langle 1\rangle)$ is isomorphic to the unit form.

Moreover, for all $i\geq m$, we have $$s(\alpha^i)=\sum_{j=0}^{m-1}\binom{i}{j}a^{i-j}.$$

\end{lem}

\begin{proof}
By Lemma \ref{unit}, it is enough to show that $s_*(\langle 1\rangle)$ is nonzero, non-degenerate, non-alternating, and that $s(x^2)$ is a square for all $x\in k[X]/(f)$.

It is easy to check that the Gram matrix $B$ of $s_*(\langle 1\rangle)$ with respect to the basis $(1,\alpha-a,\ldots,(\alpha-a)^{m-1})$ is $$B=\begin{pmatrix}
1 & 1 & \cdots & 1 & 1 \cr
1 & 1 & \cdots &1 & 0 \cr
\vdots & \vdots & \adots & \adots & \vdots \cr
\vdots & 1 &\adots & & \vdots \cr
1 & 0 & \cdots &\cdots & 0
\end{pmatrix}.$$

It follows at once that $s_*(\langle 1\rangle)$ is nonzero, non-degenerate and non-alternating. The fact that $s(x^2)$ is a square for all $x\in k[X]/(f)$ follows from Remark \ref{squares} and the fact that the diagonal entries of $B$ are squares.

Let us prove the last part of the lemma, and let $i\geq m.$ Since $(\alpha-a)^j=0$ in $k[X]/(f)$ for all $j\geq m,$ we have
 $$\alpha^i=(\alpha-a+a)^i=\sum_{j=0}^{m-1} \binom{i}{j} (\alpha-a)^j a^{i-j}.$$
 Applying the definition of $s$ yields the desired equality.
\end{proof}

\begin{prop}\label{transsepm}
Let $\pi\in k[X]$ be a monic irreducible separable polynomial of degree $d\geq 1$,  let $a\in k_s$ be a root of $\pi$, and let $L=k(a)$. Let $m\geq 1$, and let $t:L[X]/((X-a)^m)\to L$ be 
the unique $L$-linear form such that $$t(1)=t(\gamma-a)=\cdots=t((\gamma-a)^{m-1})=1,$$
where $\gamma$ is the class of $X$ in the quotient ring $L[X]/((X-a)^m)$.

Finally, let $s:k[X]/(\pi^m)\to k$ be the $k$-linear map defined by $$s(\ov{P})=\Tr_{L/k}(t(\widetilde{P})) \ \mbox{ for all }\ov{P}\in k[X]/(\pi^m),$$
where $\widetilde{Q}$ denotes the class of a polynomial $Q\in L[X]$ in $L[X]/((X-a)^m)$.

Then $s_*(\langle 1\rangle)$ is isomorphic to the unit form.

Moreover, for all $i\geq m$, we have $$s(\alpha^i)=\sum_{j=0}^{m-1}\binom{i}{j}\Tr_{L/k}(a^{i-j}).$$ 
\end{prop}

\begin{proof}
First, note that $s$ is well-defined. Indeed, if $P\in k[X]$ is a multiple of $\pi^m$, then it is a multiple of $(X-a)^m$ in $k_s[X]$, since $a$ is a root of $\pi$. In particular, $\widetilde{P}$ only depends on the class $\ov{P}$ (and not on the choice of a representative $P\in k[X]$).

By Lemma \ref{xam}, we can pick an orthonormal $L$-basis $(\widetilde{Q}_1,\ldots, \widetilde{Q}_m)$ of the $L$-vector space $L[X]/((X-a)^m)$ with respect to $t_*(\langle 1\rangle)$. By Lemma \ref{transtr}, we can also pick an orthonormal $k$-basis $(\gamma_1,\ldots,\gamma_d)$ of $L/k$ with respect to $(\Tr_{L/k})_*(\langle 1 \rangle)$, where $d=[L:k]=\deg(\pi)$. 

{\bf Claim. }For all $Q\in L[X]$, there exists $P\in k[X]$ such that $\widetilde{P}=\widetilde{Q}.$

Assume the claim is proved, and let $P_{ij}\in k[X]$ be a polynomial such that $\widetilde{P}_{ij}=\widetilde{ \gamma_i Q_j}.$
Then for all $1\leq i,r\leq d,1\leq j,s\leq m,$ we have $$s(\ov{P}_{ij}\ov{P}_{rs})=\Tr_{L/k}(t(\widetilde{P}_{ij}\widetilde{P}_{rs})) =\Tr_{L/k}(t(\widetilde{\gamma_iQ_j}\widetilde{\gamma_rQ_s})).$$
Since $t$ is $L$-linear, this yields $$s(\ov{P}_{ij}\ov{P}_{rs})=\Tr_{L/k}(\gamma_i\gamma_r t(\widetilde{Q}_j\widetilde{Q}_s))=\Tr_{L/k}(\gamma_i\gamma_r \delta_{js})=\delta_{js}\delta_{ir}.$$
It follows that the family $(\ov{P}_{ij})_{i,j}$ is orthonormal with respect to $s_*(\langle 1\rangle)$. In particular, the $P_{ij}$'s are linearly independent over $k$. Since this family has $dm$ elements, we may conclude that $(P_{ij})_{i,j}$ is an orthonormal basis of $k[X]/(\pi^m)$ with respect to $s_*(\langle 1\rangle)$.

It remains to prove the claim. Let $M/k$ be the Galois closure of $L$, let $G$ be the Galois group of $M/k$, and let $X_L$ be the set of $k$-embeddings $\sigma:L\to k_{alg}$. Note that elements of $X_L$ have their images contained in $M$, so $G$ acts by left composition on $X_L$.

 Let $Q\in L[X]$.
Since $\pi$ is separable, the polynomials $(X-\sigma(a))^m, \sigma\in X_L$ are pairwise coprime. The Chinese Remainder Theorem yields the existence of a 
unique polynomial $P\in M[X]$ of degree $<dm$ such that $P\equiv \sigma(Q) \ \mod (X-\sigma(a))^m$ for all $\sigma \in X_L$, where the congruences are viewed inside $M[X]$.

For all $\tau\in G$, we then have $\tau(P)\equiv \tau\sigma(Q) \ \mod (X-\tau\sigma(a))^m$ for all $\sigma \in X_L$.
 Since $\pi$ is irreducible, the action of $G$ on $X_L$ is transitive and we have $\tau(P)\equiv \sigma(Q) \mod (X-\sigma(a))^m$ for all $\sigma\in X_L$. Now, $\deg(\tau(P))=\deg(P)$, so by uniqueness of $P$, we get $\tau(P)=P$ for all $\tau\in G$. Hence $P\in k[X]$. By choice of $P$, we have $(X-a)^m\mid (P-Q)$ in $M[X]$. Since $P\in k[X]\subset L[X]$, and $(X-a)^m\in L[X]$, the corresponding quotient also lies in $L[X]$, and we have $P\equiv Q \mod (X-a)^m$ in $L[X]$, that is $\widetilde{P}=\widetilde{Q}$. 
 
The last part comes from Lemma \ref{xam}. This concludes the proof.
\end{proof}

\begin{lem}\label{x2nam}
Let $m\geq 2$ and $n\geq 1$. Assume that $a\in k$ is not a square, and let $f=(X^{2^n}-a)^m$. Let $s:k[X]/(f)\to k$ be the unique $k$-linear form such that 
$$s(\alpha^{2^n})=1, s(\alpha^{2^nm-1})=1, s(\alpha^j)=0 \ \mbox{ if }j\neq 2^n, 2^n m-1.$$
Then $s_*(\langle 1\rangle)$ is isomorphic to the unit form. 

Moreover, for all $i\geq 2^nm$, $s(\alpha^i)$ is equal to :

\begin{enumerate}
\item  $\ds\binom{2^nu-1}{2^nm-1}a^{u-m}$ if $i=2^nu-1$ for some $u>m$

\item  $\ds\Bigl(\sum_{j=0}^{2^n (m-1)-1}\binom{2^{n+1}u}{j}\Bigr)a^{2u}$ if $i=2^n(2u+1)$ for some $u\geq \dfrac{m-1}{2}$

\item $0$ otherwise.
\end{enumerate}
\end{lem}

\begin{proof}
Once again, we are going to use Lemma \ref{unit}. Note for the computations that $2^n<2^nm-1$ since $m\geq 2$ and $n\geq 1$.

Clearly, $s_*(\langle 1 \rangle)$ is nonzero, and non-alternating since $s(\alpha^{2n})=1$. Let us prove it is non-degenerate. Let $x\in k[X]/((X^{2^n}-a)^m)$ be such that $s(xy)=0$ for all $y\in  k[X]/((X^{2^n}-a)^m)$. Assume that $x\neq 0$. Then the ideal $(x)$ is contained in $\ker(s)$. 
Since $a$ is not a square and $k$ has characteristic two, $X^{2^n}-a$ is irreducible in $k[X]$, and the monic divisors of $(X^{2^n}-a)^m$ are the polynomials $(X^{2^n}-a)^j, \ 0\leq j\leq m$. 
Thus, the ideals of $k[X]/((X^{2^n}-a)^m)$ are the ideals $((\alpha^{2^n}-a)^j),   \ 0\leq j\leq m$. Since $(x)$ is nonzero and is not equal to the whole quotient ring (since $s$ is not identically zero), we have $(x)=((\alpha^{2^n}-a)^j)$ for some $1\leq j\leq m-1$. In particular, $(x)$ contains $(\alpha^{2^n}-a)^{m-1}$. Hence, $(\alpha^{2^n}-a)^{m-1}\in \ker(s)$ and we have $s((\alpha^{2^n}-a)^{m-1}\alpha^{2^n-1})=0$.
Now $(\alpha^{2^n}-a)^{m-1}\alpha^{2^n-1}$ is a polynomial in $\alpha$ of degree $2^nm-1$, whose coefficient of $\alpha^{2^n}$ is zero and whose coefficient of $\alpha^{2^nm-1}$ is $1$. Therefore,  $s((\alpha^{2^n}-a)^{m-1}\alpha^{2^n-1})=1$ and we get a contradiction.

It remains to prove that $s(\alpha^{2i})$ is a square for all $0\leq i\leq  2^nm -1$. This is clear if $0\leq 2i\leq 2^nm-2$.
For the other cases, we are going to prove the formulas giving $s(\alpha^i)$ for all $i\geq 2^nm$. 
Let $i\geq 2^nm,$ and let $\beta\in k_{alg}$ satisfying $\beta^{2^n}-a=0$. Note that $f=(X-\beta)^{2^nm}$.

In $k_{alg}[X]$, we have 
$$X^i=(X-\beta+\beta)^i\equiv \sum_{j=0}^{2^nm-1}\binom{i}{j}(X-\beta)^j \beta^{i-j} \ \mod (X-\beta)^{2^nm}.$$
Since the degree of $R_i=\ds\sum_{j=0}^{2^nm-1}\binom{i}{j}(X-\beta)^j \beta^{i-j}$ is $\leq 2^n m-1$, $R_i\in k_{alg}[X]$ is the remainder of the long division of $X^i$ by $(X-\beta)^{2^nm}=f$. Since $X^i$ and $f$ lie in $k[X]$, so is $R_i$. If $r_j^{(i)}\in k$ is the coefficient of $X^j$ in $R_i$, we then have $s(\alpha^i)=r_{2^n}^{(i)}+r_{2^nm-1}^{(i)}$.

Since $X^{2^n}-a$ is irreducible, $\beta$ has degree $2^n$. It follows that $1,\beta,\ldots,\beta^{2^n -1}$ are $k$-linearly independent. If $j\geq 0,$ a Euclidean division of $j$ by $2^n$ then shows that $\beta^j$ lies in $k$ if and only if $j$ is a multiple of $2^n$.

We clearly have $r_{2^nm-1}^{(i)}=\ds\binom{i}{2^nm-1}\beta^{i-(2^nm-1)}$. Note that $\ds\binom{i}{2^nm-1}=0$ or $1$ in $k.$ In particular, if  $r_{2^nm-1}^{(i)}$ is nonzero, it is equal to $\beta^{i-(2^nm-1)}$. Since  $r_{2^nm-1}^{(i)}\in k,$ we necessarily have $i\equiv - 1 \ [2^n]$. This implies that $r_{2^nm-1}^{(i)}=0$ when $i$ is even. Now, if $i=2^n u-1$ (where $u>m$), we then get $r_{2^nm-1}^{(i)}=\ds\binom{2^nu-1}{2^nm-1}a^{u-m}$.

We also have $$r_{2^n}^{(i)}=\sum_{j=2^n}^{2^n m-1}\binom{i}{j}\binom{j}{2^n}\beta^{j-2^n}\beta^{i-j}=\varepsilon_i\beta^{i-2^n},$$ where $\varepsilon_i=\ds\sum_{j=2^n}^{2^n m-1}\binom{i}{j}\binom{j}{2^n}=\binom{i}{2^n}\Bigl(\sum_{j=2^n}^{2^n m-1}\binom{i-2^n}{j-2^n}\Bigr)$.
As before, $\varepsilon_i=0$ or $1$ in $k$.
If $r_{2^n}^{(i)}\neq 0$, we then have $r_{2^n}^{(i)}=\beta^{i-2^n}\in k$, which implies that $i\equiv 0 \ [2^n]$.

Hence, if $r_{2^n}^{(i)}\neq 0$, there exists $\ell \geq m$ such that $i=2^n \ell$. In particular, $r_{2^n}^{(i)}= 0$ if $i$ is odd. Now it is well-known that $\ds\binom{2^n\ell}{2^n}$ and $\ell$ have same parity. Since $r_{2^n}^{(i)}\neq 0$, so is $\varepsilon_i$, and $\ell$ is necessarily odd. Thus $\ell=2u+1$ for some $u\geq \dfrac{m-1}{2}$.

Consequently, if $i=2^n(2u+1)$, we have $r_{2^n}^{(i)}=\varepsilon_i \beta^{2^{n+1}u}=\varepsilon_i a^{2u}$, where $$\varepsilon_i=\sum_{j=2^n}^{2^n m-1}\binom{i-2^n}{j-2^n}=\sum_{j=0}^{2^n (m-1)-1}\binom{2^{n+1}u}{j}.$$

All the desired  formulas follow at once. In particular, $s(\alpha^{2i})$ is a square for all $i\geq 0$.
This concludes the proof.
\end{proof}

\begin{prop}\label{transinsepm}
Let $\rho\in k[X]$ be a monic irreducible inseparable polynomial. Write $\rho=\pi(X^{2^n})$, where $n\geq 1$ and 
 $\pi\in k[X]$  is a monic irreducible separable polynomial of degree $d\geq 1$.
Let $a\in k_s$ be a root of $\pi$, and let $L=k(a)$. Let $t:L[X]/((X^{2^n}-a)^m)\to L$ be 
the unique $L$-linear form such that
$$t(\gamma^{2^n})=1, t(\gamma^{2^nm-1})=1, t(\gamma^j)=0 \ \mbox{ if }j\neq 2^n, 2^n m-1.$$
where $\gamma$ is the class of $X$ in the quotient ring $L[X]/((X^{2^n}-a)^m)$.

Finally, let $s:k[X]/(\pi^m)\to k$ be the $k$-linear map defined by $$s(\ov{P})=\Tr_{L/k}(t(\widetilde{P})) \ \mbox{ for all }\ov{P}\in k[X]/(\pi^m),$$
where $\widetilde{Q}$ denotes the class of a polynomial $Q\in L[X]$ in $L[X]/((X^{2^n}-a)^m)$.

Then $s_*(\langle 1\rangle)$ is isomorphic to the unit form.

Moreover, for all $i\geq 2^nm$, $s(\alpha^i)$ is equal to :

\begin{enumerate}
\item  $\ds\binom{2^nu-1}{2^nm-1}\Tr_{L/k}(a^{u-m})$ if $i=2^nu-1$ for some $u> m$

\item  $\ds\Bigl(\sum_{j=0}^{2^n (m-1)-1}\binom{2^{n+1}u}{j}\Bigr)\Tr_{L/k}(a^{2u})$ if $i=2^n(2u+1)$ for some $u\geq \dfrac{m-1}{2}$

\item $0$ otherwise.
\end{enumerate}
\end{prop}

\begin{proof}
If $\sigma_1,\sigma_2\in X_L$ are two distinct $k$-embeddings of $L$ into $k_{alg}$, the polynomials $X^{2^n}-\sigma_1(a)$ and $X^{2^n}-\sigma_2(a)$ are coprime (since $\pi$ is separable), and thus so are $(X^{2^n}-\sigma_1(a))^m$ and $(X^{2^n}-\sigma_2(a))^m$.
Reasoning as in the proof of Proposition \ref{transsepm}, we may show that, given $Q\in L[X]$, there exists $P\in k[X]$ such that $\widetilde{P}=\widetilde{Q}$. The rest of the proof is then identical to the proof of this very same proposition.
\end{proof}

\section{Proof of the main results and examples}

We now come to the proofs of the main results of this paper, that we state again for the convenience of the reader.

\begin{thm}\label{main1}
Let $f\in k[X]$ be a monic polynomial of degree $n\geq 1$. Then  $f$ is the minimal polynomial of a symmetric matrix with entries in $k$ if and only if $f$ is not the product of pairwise distinct monic irreducible inseparable polynomials.

In this case, $f$ is the minimal polynomial of a symmetric matrix of $\M_n(k)$.

In particular, if $k$ is perfect, any monic  polynomial $f\in k[X]$ of degree $n\geq 1$  is the minimal polynomial of a symmetric  matrix of $\M_n(k)$.
\end{thm} 

\begin{proof}
We will make use of the following easy fact : 
if $E_1,\ldots,E_s$ are finite dimensional commutative $k$-algebras, we have $$\mu_x={\rm lcm}(\mu_{x_1},\ldots,\mu_{x_s}) \ \mbox{ for all }x=(x_1,\ldots,x_s)\in E_1\times\cdots\times E_s.$$

Assume first that $f$ has at least one separable irreducible divisor, and write
$f=\pi_1^{m_1}\cdots \pi_r^{m_r}g,$ where $r\geq 1, m_1,\ldots,m_r\geq 1$, $\pi_1,\ldots,\pi_r$ are monic separable irreducible polynomials, and $g$ is the product of the inseperable irreducible divisors of $f$ (with multiplicities). Note that $g$ lies in $k[X^2]$ and is coprime to each $\pi_i^{m_i}$.
By Proposition \ref{transsepm}, for all $1\leq i \leq r$,  there exists a $k$-linear form $s_i:k[X]/(\pi_i^{m_i})\to k$ such that $(s_i)_*(\langle 1\rangle)$ is isomorphic to the unit form. Since $g\in k[X^2],$ by Lemma \ref{even}, there exists a $k$-linear form $s':k[X]/(g)\to k$ such that $s'_*(\langle 1\rangle)$ is hyperbolic. By Lemma \ref{sumtrans}, there exists a $k$-linear map $s:k[X]/(\pi_1^{m_1})\times\cdots\times k[X]/(\pi_r^{m_r})\times k[X]/(g)\to k$ such that $s_*(\langle 1\rangle)$ is isomorphic to the unit form. 
If $\alpha_i$ is the class of $X$ in $k[X]/(\pi_i^{m_i})$ and $\alpha'$ is the class of $X$ in $k[X]/(g))$, we have $\mu_{\alpha_i}=\pi_i^{m_i}$ and $\mu_{\alpha'}=g$. Let $x=(\alpha_1,\ldots,\alpha_r,\alpha')$. Since $\pi_1^{m_1},\ldots,\pi_r^{m_r}$ and $g$ are pairwise coprime, we have $$\mu_x={\rm lcm}(\pi_1^{m_1},\ldots,\pi_r^{m_r},g)=\pi_1^{m_1}\cdots \pi_r^{m_r}g=f.$$
By Lemma \ref{kxftransfer}, $f$ is the minimal polynomial of a symmetric matrix of $\M_n(k)$.

Assume now that all irreducible divisors of $f$ are inseparable. If $f$ is the product of pairwise distinct monic irreducible inseparable polynomials, then $f$ is not the minimal polynomial of a symmetric matrix of any size by Proposition \ref{prodrho}.

Otherwise, we may write
$f=\rho^m g$, where $m\geq 2$, $\rho$ is a monic irreducible inseparable polynomial, and $g$ is the product of the other  inseperable irreducible divisors of $f$ (with multiplicities). Once again, $g$ lies in $k[X^2]$ and is coprime to $\rho^m$. By Proposition \ref{transinsepm}, there exists a $k$-linear form $s:k[X]/(\rho^m)\to k$ such that $s_*(\langle 1\rangle)$ is isomorphic to the unit form. By Lemma \ref{even}, there exists a $k$-linear form $s':k[X]/(g)\to k$ such that $s'_*(\langle 1\rangle)$ is hyperbolic.  Now, we may finish the argument as before, and this concludes the proof.
\end{proof}

\begin{coro}\label{coromain1}
Let $f\in k[X]$ be a monic  polynomial of degree $n\geq 1$. Then $f$ is the characteristic polynomial of a symmetrix matric of $\M_n(k)$ if and only if is not the product of pairwise distinct monic irreducible inseparable polynomials.

In particular, if $k$ is perfect, any monic polynomial $f\in k[X]$ of degree $n\geq 1$ is the characteristic polynomial of a symmetric matrix of $\M_n(k)$.
\end{coro}

\begin{proof}
Assume that $f$ is not the product of pairwise distinct monic irreducible inseparable polynomials. Theorem \ref{main1} gives the existence of a symmetric matrix $M\in \M_n(k)$ such that $\mu_M=f$. Since $\mu_M$ and $\chi_M$ are both monic polynomials of degree $n$ and $\mu_M\mid \chi_M$, we have $\chi_M=\mu_M=f$.

Assume now that $f$ is  the product of pairwise distinct monic irreducible inseparable polynomials, and suppose that there exists a symmetric matrix $M\in \M_n(k)$ such that $\chi_M=f$. Since $\mu_M$ and $\chi_M$ have the same irreducible divisors, the hypothesis on $f$ implies that $\mu_M=\chi_M=f$,  contradicting the previous theorem.
\end{proof}

\begin{thm}\label{main2}
Let $\alpha\in k_{alg}$ be an algebraic element of degree $n$, with minimal polynomial $f.$ Then :

\begin{enumerate}
 \item if $f$ is separable, $\alpha$ is the eigenvalue of a symmetric matrix of $\M_n(k)$;
 
 \item if $f$ is inseparable, $\alpha$ is the eigenvalue of a symmetric matrix of $\M_{n+1}(k)$, but not of any symmetric matrix of $\M_n(k)$.
\end{enumerate}

In particular, if $k$ is perfect, any algebraic element of degree $n$ is the eigenvalue of a symmetric  matrix of $\M_n(k)$.
\end{thm}

\begin{proof}
If $f$ is separable, we may apply the previous corollary to get a symmetric matrix $M\in\M_n(k)$ such that $\chi_M=f$. In particular, $\alpha$ is an eigenvalue of $M$.

Assume now that $f$ is inseparable, and suppose that there exists a symmetric matrix $M\in \M_n(k)$ such that $\alpha$ is an eigenvalue of $M$. Since $\chi_M(\alpha)=0$ and $\chi_M\in k[X]$, we get $f\mid \chi_M$. For degree reasons, we have $\chi_M=f$, contradicting Corollary \ref{coromain1}.

Since $f\in k[X^2]$ in this case, by Lemma \ref{even},  there exists a $k$-linear form $s':k[X]/(f)\to k$ such that $s'_*(\langle 1\rangle)$ is hyperbolic.
Note now that the $k$-linear form $s_0:k[X]/(X)\to k$ induced by evaluation at $0$ obviously satifies $(s_0)_*(\langle 1\rangle)\simeq \langle 1\rangle$.
We may now conclude as in the proof of Theorem \ref{main1} that there exists a symmetric matrix $M\in \M_{n+1}(k)$ such that $\mu_M=Xf$. Degree considerations show that $\chi_M=Xf$. In particular, $\alpha$ is an eigenvalue of $M$.
This concludes the proof.
\end{proof}

\begin{rem}
In all the previous results, the symmetric matrices $M$ which we have constructed in the proofs are cyclic, that is, they satisfy $\chi_M=\mu_M$.
\end{rem}

We now summarize the method to compute explicitely a symmetric matrix of $\M_n(k)$ of given minimal polynomial $f$ of degree $n\geq 1$, provided that we know the factorization of $f$. Of course, we assume that $f$ is not the product of pairwise distinct inseparable irreducible polynomials.

We have to consider two cases :

\begin{enumerate}

 \item If $f$ has at least one separable irreducible divisor, write $f=\pi_1^{m_1}\cdots \pi_r^{m_r}g,$ where $r\geq 1, m_1,\ldots,m_r\geq 1$, $\pi_1,\ldots,\pi_r$ are monic separable irreducible polynomials, and $g$ is a polynomial of $k[X^2]$ coprime to each $\pi_i$.

Note that we can only consider the separable irreducible divisors $\pi_1,\ldots,\pi_r$ such that $m_1,\ldots,m_r$ are odd, and include the other ones in the polynomial $g$.

For all $(x_1,\ldots,x_r,x')\in E=k[X]/(\pi_1^{m_1})\times \cdots\times  k[X]/(\pi_r^{m_r})\times k[X]/(g)$, set $$s(x_1,\ldots,x_r,x')=s_1(x_1)+\cdots+s_r(x_r)+s'(x'),$$
where $s_i:k[X]/(\pi_i^{m_i})\to k$ is the $k$-linear map defined in Proposition \ref{transsepm} and $s':k[X]/(g)\to k$ is the $k$-linear map defined in Lemma \ref{even}.

In this case, we will denote by $\mathcal{E}$ the $k$-basis $$(1,0,\ldots,0),\ldots,(\alpha_1^{m_1d_1-1},0,\ldots,0),\ldots,(0,\ldots,0,1),\ldots,(0,\ldots,0,\alpha^{\prime d'-1}),$$
where $\alpha_i$ is the class of $X$ in $k[X]/(\pi_i^{m_i})$, $\alpha'$ is the class of $X$ in $k[X]/(g)$, $d_i=\deg(\pi_i)$ and $d'=\deg(g)$.

We also set $C=\begin{pmatrix}
C_{\pi_1^{m_1}} & & & \cr
                & \ddots &  & \cr
                &        & C_{\pi_r^{m_r}} & \cr
                 &        &                 & C_g
\end{pmatrix}$, where $C_h$ is the companion matrix of a monic polynomial $h\in k[X]$.

\item If all irreducible divisors of $f$ are inseparable, write $f=\rho^m g$, where $\rho$ is irreducible and inseparable, $m\geq 2$, and $g$ is a polynomial of $k[X^2]$ coprime to 
$\rho$. 

For all $(x,x')\in E=k[X]/(\rho^m)\times k[X]/(g)$, set $$s(x,x')=s(x)+s'(x'),$$ where $s:k[X]/(\rho^m)\to k$  is the $k$-linear map defined in Proposition \ref{transinsepm} and $s':k[X]/(g)\to k$ is the $k$-linear map defined in Lemma \ref{even}.

In this case, we will denote by $\mathcal{E}$ the $k$-basis $$(1,0),\ldots,(\alpha^{md-1},0),\ldots,(0,1),\ldots,(0,\alpha^{\prime d'-1}),$$
where $\alpha$ is the class of $X$ in $k[X]/(\rho^m)$, $\alpha'$ is the class of $X$ in $k[X]/(g)$, $d=\deg(\rho)$ and $d'=\deg(g)$.

We also set $C=\begin{pmatrix}
C_{\rho^m} & \cr & C_g
\end{pmatrix}$.
\end{enumerate}

Then, we may compute an orthonormal basis $\mathcal{B}$ of $E$ with respect to $s_*(\langle 1 \rangle),$ for example starting from the expression of $s(xy)$ in terms of the coordinates of $x$ and $y$ in the basis $\mathcal{E}$, and applying the algorithm explained in Section \ref{secbil}. In case (a), to avoid working with too many variables, one may first compute an orthonormal basis for each of the bilinear spaces $(k[X]/(\pi_1^{m_i}), (s_i)_*(\langle 1\rangle)),$  ($1\leq i\leq r-1$) and  $(k[X]/(\pi_r^{m_r}) \times k[X]/(g), (s_r\oplus s')_*(\langle 1\rangle)),$
and glue them into an orthonormal basis $\mathcal{B}$ of $E$.

By choice of $\mathcal{E}$, the representative matrix of left multiplication by $(\alpha_1,\ldots,\alpha_r,\alpha')$ in case (a), or by $(\alpha, \alpha')$ in case (b), is the matrix $C$. Hence, the matrix representation of this endomorphism with respect to the orthonormal basis $\mathcal{B}$ is $M=P^{-1}CP=QCQ^{-1}$ (with the notation of the algorithm explained in Section \ref{secbil}).

The matrix $M$ is then symmetric, with minimal polynomial (and characteristic polynomial) equal to $f$.

\begin{ex}
Let $k=\mathbb{F}_2(t)$, and let $f=(X^2+X+t)^3$. We will use the $k$-linear map defined in Proposition \ref{transsepm}. The representative matrix of $s_*(\langle 1 \rangle)$ with respect to the basis $(1,\alpha,\ldots,\alpha^5)$ is $$S=\begin{pmatrix}
     0&       1 &      1&       t+1 &      1 &  t^2+t\cr
      1     &  1&       t+1 &      1&   t^2+t&       t^2\cr
      1 &      t+1&       1 &  t^2 + t &      t^2 &t^3 + 1\cr
      t+1   &    1  & t^2 + t  &     t^2& t^3 + 1 &    1\cr
    1 &  t^2 + t &      t^2 &t^3 + 1 &    1& t^4 + t^2+t\cr
  t^2 + t    &   t^2& t^3 + 1    & 1& t^4 + t^2+t       &t^4+t^2
\end{pmatrix}.$$
If $x=\ds\sum_{i=0}^5x_i\cdot \alpha^i$ and $y=\ds\sum_{i=0}^5y_i\cdot \alpha^i$, then $s(xy)$ contains the terms
$$x_1y_1+x_1(y_0+(t+1)y_2+y_3+(t^2+t)y_4+t^2y_5)+y_1(x_0+(t+1)x_2+x_3+(t^2+t)x_4+t^2x_5),$$
which may be rewritten as $$(x_0+ x_1+(t+1)x_2+x_3+(t^2+t)x_4+t^2x_5)(y_0+y_1+(t+1)y_2+y_3+(t^2+t)y_4+t^2y_5)$$ $$+
(x_0+(t+1)x_2+x_3+(t^2+t)x_4+t^2x_5)(y_0+(t+1)y_2+y_3+(t^2+t)y_4+t^2y_5).$$
We then have $s_*(\langle 1\rangle)=\varphi_1\bullet\varphi_1+b_1,$ where $$\nfuncp{\varphi_1}{E}{k}{x}{x_0+x_1+(t+1)x_2+x_3+(t^2+t)x_4+t^2x_5}$$
and $b_1$ is a symmetric bilinear form that we will not write here.
After simplification, one can see that $b_1(x,y)$ contains the terms $$x_0y_0 +x_0(ty_2+ty_3+(t^2+t+1)y_4+ty_5)+y_0(tx_2+tx_3+(t^2+t+1)x_4+tx_5),$$
that is $$(x_0+tx_2+tx_3+(t^2+t+1)x_4+tx_5)(y_0+ty_2+ty_3+(t^2+t+1)y_4+ty_5)$$ $$+(tx_2+tx_3+(t^2+t+1)x_4+tx_5)(ty_2+ty_3+(t^2+t+1)y_4+ty_5).$$
Hence,  $s_*(\langle 1\rangle)=\varphi_1\bullet\varphi_1+\varphi_2\bullet\varphi_2+b_2,$ where $$\nfuncv{\varphi_2}{E}{k}{x}{x_0+tx_2+tx_3+(t^2+t+1)x_4+tx_5}$$ and $$b_2(x,y)=x_2y_3+x_2y_5+x_3y_2+x_3y_3+x_3y_4+x_3y_5+x_4y_3+x_5y_2+x_5y_3.$$
Now, $b_2(x,y)$ contains the terms $$x_3y_3+ x_3(y_2+y_4+y_5)+y_3(x_2+x_4+x_5),$$
so  $s_*(\langle 1\rangle)=\varphi_1\bullet\varphi_1+\varphi_2\bullet\varphi_2+\varphi_3\bullet \varphi_3+b_3,$
where $$\nfuncv{\varphi_3}{E}{k}{x}{x_2+x_3+x_4+x_5}$$ and $$b_3(x,y)=x_2y_2+x_2y_4+x_4y_2+x_4y_4+x_4y_5+x_5y_4+x_5y_5.$$
Since $b_3(x,y)$ contains $x_2y_2+x_2y_4+x_4y_2$, we may write  $s_*(\langle 1\rangle)=\varphi_1\bullet\varphi_1+\varphi_2\bullet\varphi_2+\varphi_3\bullet \varphi_3+\varphi_4\bullet\varphi_4+b_4,$ where $$\nfuncv{\varphi_4}{E}{k}{x}{x_2+x_4}$$ and $$b_4(x,y)=x_4y_5+x_5y_4+x_5y_5=(x_4+x_5)(y_4+y_5)+x_4y_4.$$ Finally, $s_*(\langle 1\rangle)=\ds\sum_{i=1}^6\varphi_i\bullet\varphi_i$, where  $$\nfuncv{\varphi_5}{E}{k}{x}{x_4+x_5}$$
and $$\nfuncp{\varphi_6}{E}{k}{x}{x_4}$$

The corresponding matrix $Q$ is then $$Q=\begin{pmatrix}
1 & 1& t+1&1 & t^2+t & t^2\cr
1 & 0 & t & t & t^2+t+1 & t\cr
0 & 0 & 1 &1 &1 &1 \cr
0&0&1&0&1&0\cr
0&0&0&0&1&1\cr
0&0&0&0&1&0
\end{pmatrix}.$$
If $C$ is the companion matrix of $f$, we then get $$M=QCQ^{-1}=\begin{pmatrix}
t+1& t& 1& 1& 0& 0\cr
t& t& 1& 1& 0& 0\cr
1& 1& t& t& 1& 1\cr
1& 1& t& t+1& 1& 1\cr
0& 0& 1& 1& t+1& t\cr
0& 0& 1& 1& t& t
\end{pmatrix}.$$
\end{ex}

\begin{ex}
Let $k=\mathbb{F}_2(t)$, and let $f=(X^2+t)^3$. We will use the $k$-linear map defined in Lemma \ref{x2nam}. The representative matrix of $s_*(\langle 1 \rangle)$ with respect to the basis $(1,\alpha,\ldots,\alpha^5)$ is $$S=\begin{pmatrix}
     0&       0 &      1&       0 &      0 &  1\cr
      0     &  1&       0 &      0&  1&       t^2\cr
      1 &      0&       0 &  1&      t^2 &t\cr
      0   &    0  & 1  &     t^2& t &    0\cr
    0 &  1 &      t^2 &t&    0& 0\cr
  1    &   t^2& t    & 0& 0       &t^4
\end{pmatrix}.$$

This time, we may prove that $s_*(\langle 1 \rangle)=\ds\sum_{i=1}^4\varphi_i\bullet\varphi_i+b_4,$ where 
$$\nfuncv{\varphi_1}{E}{k}{x}{x_1+x_4+t^2x_5} \ \nfuncv{\varphi_2}{E}{k}{x}{t^2x_2+tx_3+x_4+t^2x_5}$$ $$ \nfuncv{\varphi_3}{E}{k}{x}{t^{-2}x_0+t^2 x_2+(t+t^{-2})x_3+(t^2+t^{-1})x_5}$$ 
$$\nfuncv{\varphi_4}{E}{k}{x}{t^{-2}x_0+(t+t^{-2})x_3+t^{-1}x_5}$$
and $b_4(x,y)=x_3y_5+y_3x_5.$

Since $b_4$ is alternating, we replace $\varphi_4$ by $$\nfuncv{\varphi_4}{E}{k}{x}{t^{-2}x_0+(t+t^{-2}+1)x_3+t^{-1}x_5}$$
and we set  $$\nfuncv{\varphi_5}{E}{k}{x}{t^{-2}x_0+(t+t^{-2})x_3+(t^{-1}+1)x_5}$$
and  $$\nfuncp{\varphi_6}{E}{k}{x}{t^{-2}x_0+(t+t^{-2}+1)x_3+(t^{-1}+1)x_5}$$

The corresponding matrix $Q$ is $$Q=\begin{pmatrix}
0 & 1 & 0 & 0 & 1 & t^2 \cr
0 & 0 & t^2 & t & 1 & t^2 \cr
t^{-2} & 0 & t^2 & (t+t^{-2}) & 0 & t^2+t^{-1} \cr
t^{-2} & 0 & 0 & t+t^{-2}+1 & 0 & t^{-1} \cr
t^{-2} & 0 & 0 & t+t^{-2} & 0 & t^{-1}+1 \cr
t^{-2} & 0 & 0 & t+t^{-2}+1 & 0 & t^{-1}+1 \cr
\end{pmatrix}.$$

If $C$ is the companion matrix of $f$, we have $M=QCQ^{-1}$, and we finally get 

$$M=t^{-4}\begin{pmatrix}
0& t^6& t^6& 0& 0& 0\cr
t^6& 0& t^3& t^3& t^4+t^3& t^4+t^3\cr
t^6& t^3& 1& t^2+1& t^4+1& t^4+t^2+1\cr
0& t^3& t^2+1& t^4+1& t^4+t^2+1& 1\cr
0& t^4+t^3& t^4+1& t^4+t^2+1& t^5+1& t^5+t^2+1\cr
0& t^4+t^3& t^4+t^2+1& 1& t^5+t^2+1& t^5+t^4+1
\end{pmatrix}.$$

\end{ex}

Note that the proofs of Propositions \ref{transsepm} and \ref{transinsepm} provide a way to compute an orthonormal basis of $k[X]/(h^m)$ for the appropriate $k$-linear form  when $h$ is an irreducible  (separable or inseparable) polynomial (with $m\geq 2$ if $h$ is inseparable), which reduces the number of variables to be manipulated at the same time :

\begin{enumerate}
\item If $h=\pi(X^{2^n})$ (where $n=0$ if $h$ is separable), and if $L=k(a)$, where $a\in k_s$ is a fixed root of $\pi$, compute an orthonormal basis $(\gamma_1,\ldots,\gamma_d)$ of $L$ with respect to $(\Tr_{L/k})_*(\langle 1 \rangle),$ using for example the method explained in Section \ref{secbil}.

\item Compute an orthonormal basis $(\widetilde{Q}_1,\ldots \widetilde{Q}_{2^nm})$ of $L[X]/((X^{2^n}-a)^m)$  for $t_*(\langle 1\rangle)$, where $t:L[X]/(X^{2^n}-a)^m\to L$ is the appropriate $L$-linear form (where $n=0$ if $h$ is separable).

\item For $1\leq i\leq d$ and $1\leq j\leq 2^nm$, compute the unique polynomial $P_{ij}\in L[X]$ of degree $\leq 2^nmd-1$ satisfying $P_{ij}\equiv \gamma_i Q_j \ \mod (X^{2^n}-\sigma(a))^m$ for all $\sigma\in X_L$. This polynomial lies in fact in $k[X]$, and $(\ov{P}_{ij})_{i,j}$ is the desired orthonormal basis.
\end{enumerate}

\begin{ex}
Let $k=\mathbb{F}_2(t)$, and let $f=(X^2+X+t)^3$. Let $a\in k_s$ be a root of $\pi=X^2+X+t$, and let $L=k(a)$. It is easy to check that $\Tr_{L/k}(1)=0$ and $\Tr_{L/k}(a)=1$. In particular, $\Tr_{L/k}(a^2)=\Tr_{L/k}(a)^2=1$ and $\Tr_{L/k}(a(1+a))=\Tr_{L/k}(t)=0$.
Hence, $(a,1+a)$ is an orthonormal basis of $L$ with respect to $(\Tr_{L/k})_*(\langle 1 \rangle)$.

Now, let $\gamma=\widetilde{X}\in L[X]/((X-a)^3)$. As observed in the proof of Lemma \ref{xam}, the representative matrix of $t_*(\langle 1\rangle)$ in the $L$-basis $(1,\gamma-a, (\gamma-a)^2)$ is $$\begin{pmatrix}
1 & 1 & 1 \cr 1 & 1 & 0 \cr 1 & 0 & 0
\end{pmatrix}.$$ Hence, if $x=x_0+x_1(\gamma -a )+x_2(\gamma-a)^2$ and $y=y_0+y_1(\gamma -a )+y_2(\gamma-a)^2,$ we have $$t(xy)=x_0y_0+x_0y_1+x_0y_2+y_0x_1+y_0x_2+x_1y_1.$$ 

It is not difficult to see that  $t_*(\langle 1\rangle)=\ds\sum_{i=1}^3 \varphi_i\bullet\varphi_i,$ where $$\nfuncv{\varphi_1}{E}{k}{x}{x_0+x_1+x_2} \ \nfuncv{\varphi_1}{E}{k}{x}{x_1+x_2} \ \mbox{ and } \ \nfuncp{\varphi_1}{E}{k}{x}{x_1}$$

Inverting the matrix $\begin{pmatrix}
1 & 1 & 1 \cr 0 & 1 & 1 \cr 0 & 1 & 0
\end{pmatrix}$ yields the orthonormal basis $$(1,1+(\gamma-a)^2, (\gamma-a)+(\gamma-a)^2).$$

Set $Q_1=1, Q_2=1+(X-a)^2$ and $Q_3=(X-a)+(X-a)^2$. Using a computer algebra system, we find that the polynomials $$aQ_1,(1-a)Q_1,aQ_2, (1-a)Q_2, aQ_3, (1-a)Q_3$$ lift to the following polynomials of $k[X]$ respectively :

$$P_{1,1}= X^4+ (t^2+t), \ P_{2,1}=X^4+ (t^2+t+1), \ P_{1,2}= X^5+X^4+X^3+ t^2X+(t^2+t)$$
$$P_{2,2}=X^5+X^3+X^2+t^2X+t+1, \ P_{1,3}=X^4+X^3+tX+t^2$$ $$P_{2,3}=X^4+X^3+X^2+(t+1)X+t^2+t .$$

The matrix of the orthonormal basis $(\ov{P}_{1,1},\ov{P}_{2,1},\ldots,\ov{P}_{2,3})$ in the basis $(1,\alpha,\ldots,\alpha^5)$ is thus $$P=\begin{pmatrix}
t^2+t & t^2+t+1 & t^2+t & t+1 & t^2 & t^2+t \cr 0 & 0 & t^2 & t^2 & t & t+1 \cr 0 & 0 & 0 & 1 & 0 & 1 \cr 0 & 0 & 1 & 1 & 1 & 1 \cr 1 & 1 & 1 & 0 & 1 & 1 \cr 0 & 0 & 1 & 1 & 0 & 0
\end{pmatrix}.$$ If $C=C_f$, we get this time $$M=P^{-1}CP=\begin{pmatrix}
t & t & 1 & 0 & 1 & 0\cr t & t+1 & 0 & 1 & 0 & 1 \cr 1 & 0 & t & t & 1 & 0 \cr 0 & 1 & t & t+1 & 0 & 1 \cr 1 & 0 & 1 &  0 & t+1 & t \cr 0 & 1 & 0 & 1 & t & t 
\end{pmatrix}.$$
\end{ex}

To conclude, we would like to give a result which may be useful for computations when the polynomial $g$ in the factorization of $f$ is a square (in this case, one may reduce the computation of an orthonormal basis on $k[X]/(f)$ to the computation of an orthonormal basis on each factor).

\begin{lem}
Assume that $f=X^{2m}+a_m^2 X^{2m-2}+\cdots +a_1^2X^2+a_0^2$, where $a_0\neq 0$, and let $s:k[X]/(f)\to k$ be the unique $k$-linear form such that $$s(1)=1, s(\alpha)=\cdots=s(\alpha^{2m-1})=0.$$
Then $s_*(\langle 1\rangle)$ is isomorphic to the unit form.
\end{lem}

\begin{proof}
By Lemma \ref{unit}, it is enough to show that $s_*(\langle 1\rangle)$ is nonzero, non-degenerate, non-alternating, and that $s(x^2)$ is a square for all $x\in k[X]/(f)$.

Since $s(1)=1$, $s_*(\langle 1 \rangle)$ is nonzero, and non-alternating. Let us shows that it is non-degenerate. 

Let $x=\ds\sum_{i=0}^{2m-1}\lambda_i\cdot\alpha^i\in \ker(s_*(\langle 1\rangle)$. Then $s(x 1)=s(x)=\lambda_0=0.$ Assume we proved that $\lambda_0=\cdots=\lambda_j=0$ for some $0\leq j\leq 2m-2$. Since $a_0\neq 0$, $\alpha$ is invertible in the quotient ring $k[X]/(f)$ (since $X$ is coprime to $f$).

We then have $0=s(x\alpha^{-j-1})=s(\ds\sum_{i=j+1}^{2m-1}\lambda_i\cdot \alpha^{i-j-1})=s(\ds\sum_{i=0}^{2m-j-2}\lambda_{i+j+1}\alpha^i)=\lambda_{j+1}$.
It follows by induction that all the $\lambda_i'$s are zero, that is $x=0$.

Let us finally prove that $s(x^2)$ is a square for all $x\in k[X]/(f)$.  For $i\geq 0$, let $R$ be  the remainder of the long division of $X^{2i}$ by $f$, and note that $X^{2i}$ and $f$ are squares. We claim that $R$ is a square. Indeed, since $K$ has characteristic two, $K^2$ is a field  and $(K[X])^2=K^2[X^2]$, so we may just consider $X^2$ as the variable  and $K^2$ as the new field of coefficients to conclude.

If $R=S^2$, then we have $\deg(S)\leq m-1$. Writing $S=\ds\sum_{j=0}^{m-1}s_j X^j,$ we get  $\alpha^{2i}=R(\alpha)=\ds\sum_{j=0}^{m-1}s_j^2 \alpha^{2j}$, and thus $s(\alpha^{2i})=s_0^2$ is a square. By Remark \ref{squares}, we get the desired result.
\end{proof}

\bigskip

{\it Acknowledgements. }We would like to thank warmly the anymomous referee for his/her careful reading of a previous version of this paper and his/her insightful remarks and suggestions.

\end{document}